\documentstyle[11pt]{article}

\normalbaselines
\def\section#1{\addtocounter{section}{1}\begin{center}\small\thesection. #1\end{center}}

\begin{document}

\newtheorem{lema}{Lemma}[section]
\newtheorem{prop}[lema]{Proposition}
\newtheorem{teor}[lema]{Theorem}

\font\ddpp=msbm10 %scaled \magstep 1 %Caracteres "doble palo".
\def\R{\hbox{\ddpp R}}    % numeros reales
\def\C{\hbox{\ddpp C}}     % numeros complejos
\def\L{\hbox{\ddpp L}}    % Lorentz-Minkowski
\def\Z{\hbox{\ddpp Z}}    %números enteros
\def\S{\hbox{\ddpp S}}    %de Sitter

\newcommand{\D}{{\cal D}}
\newcommand{\pD}{{\partial \cal D}}
\newcommand{\PP}{{\cal S}}
\newcommand{\CC}{{\cal C}}
\newcommand{\NN}{{\cal N}}
\newcommand{\KK}{{\cal K}}

\newcommand{\M}{{\cal M}}
\newcommand{\be}{\begin{equation}}
\newcommand{\ee}{\end{equation}}
\newcommand{\N}{{\bf N}}
\newcommand{\dem}{{\bf Proof: }}
\newcommand{\inte}{\int_{0}^{1}}
\newcommand{\eps}{\epsilon}
\newcommand{\<}{\langle}
\renewcommand{\>}{\rangle}
\renewcommand{\(}{\left(}
\renewcommand{\)}{\right)}

% \newfont{\n}{msym10}       

%\bibliographystyle{unsrt} 

\vbox {\vspace{20mm}} %Leave space at the top on the first page.

\begin{center}
{\xivpt GEODESIC CONNECTEDNESS OF SEMI-RIEMANNIAN MANIFOLDS.}  \\[2mm]
 MIGUEL S\'ANCHEZ \\[1mm] 
\baselineskip=4mm
{\ixpt Departamento de Geometr\'{\i}a y Topolog\'{\i}a, Facultad de Ciencias,
Universidad de Granada, \\ Avda. Fuentenueva s/n, 18071-Granada, Spain} 
\end{center}

\hspace{2mm}{\parbox{150mm}{\ixpt \baselineskip=4mm {\it Key words and phrases:}  Riemannian, Lorentzian and affine manifolds, 
geodesic connectedness, convex boundary, pseudoconvexity, space of geodesics, Lorentzian and affine torus, spaceform, variational methods, critical points, Ljusternik-Schnirelman theory, stationary and splitting manifolds, multiwarped and Generalized Robertson-Walker spacetimes, Brower's topological degree.}}

\def\thefootnote{}
\footnotetext{\hspace*{-4.5mm}{\ixpt  The author thanks Prof. R. Bartolo, J.L. Flores and P.E. Ehrlich for a careful reading of the manuscript.
This work has been partially supported by the DGICYT Grant PB-97-0784-C03-01.}}

\begin{center}
\section{INTRODUCTION AND RIEMANNIAN MANIFOLDS}
\end{center}
\setcounter{equation}{0}

\noindent The problem of geodesic connectedness in semi-Riemannian manifolds (i.e. the question whether each two points of the manifold can be joined by a geodesic) has been widely studied from very different viewpoints. Our purpose is to review these semi-Riemannian techniques, and possible extensions. 
In the Riemannian case, 
it is natural to state this problem on (incomplete) manifolds with (possibly non-smooth) boundary, and we will discuss different  conditions on this boundary in the remainder of  this Section. In this case, geometrical as well as variational methods are appliable, and accurate results can be obtained by using the associated distance and related properties of  positive-definiteness.  
For Lorentzian manifolds, the results cannot be so general, and very different techniques have been introduced which are satisfactory for some particular classes of Lorentzian manifolds. We start by considering several geometrical notions appliable to affine manifolds and, thus, to all semi-Riemannian manifolds, Section~2. Recall that variational methods are not appliable to affine manifolds, at least in a standard way: geodesics are the critical points of the action functional, but  a metric tensor must be provided for the definition of this functional. 
In Section 3 some general facts on geodesic connectedness of Lorentzian manifolds are pointed out, and  classical results about  connectedness of spaceforms, which rely in the properties of  actions of isommetry groups, are summarized.  
In Section 4 we discuss variational methods applied to Lorentzian manifolds, which have been shown to be useful, mainly, to study stationary and splitting manifolds, with or without boundary. Finally, in Section 5 recent results, based on topological arguments and appliable to multiwarped spacetimes, are explained.  

\vspace{2mm}

\noindent As pointed out by Gordon \cite{Go}, the problem of the geodesic connectedness of a Riemannian manifold is important not only in its own right but also because of its relation, via the Jacobi metric, with the problem of connecting the points   of the manifold by means of the trajectories determined by an autonomous potential. Moreover, it is also important because   stronger results on geodesic connectedness among all semi-Riemannian manifolds are found for  definite metrics. In fact, the following concept is especially useful: a Riemannian metric will be said to be {\it convex} if each two of its points can be joined by means of a distance-minimizing geodesic (not necessarily unique). By the Hopf-Rinow theorem all complete Riemannian manifolds are convex, and we  discuss now when an incomplete one is either convex or geodesically connected.

\vspace{2mm}

\noindent We start with the simplest case. Let   $(\M,\< \cdot, \cdot \>)$ be a $n$-dimensional Riemannian manifold (all the manifolds will be assumed to be connected and smooth, even though at most differentiability $C^4$ will be needed), and $\D \subset \M$ an open (connected) domain with differentiable (smooth) boundary $\pD$; put $\widetilde \D = \D \cup \pD$.  We will {\it not} assume a priori that $\M$ is complete, because this assumption is not especially relevant for the indefinite case. The following definitions become natural:

(I)  $\partial \D$ is {\it infinitesimally   convex at}
$p\in \partial \D$ (IC$_p$) if the second fundamental form $\sigma_p$, with respect
to the interior normal, is positive semidefinite. When $\sigma_p$ is positive definite we say that $\partial \D$ is {\it strictly locally   convex at}
$p$ (SIC$_p$). Equivalently, take any  differentiable function 
$\phi: U \cap \widetilde \D \longrightarrow\R$, where
$ U\subset \M$ is a neighborhood of $p$
 such that

%%\vspace*{-3mm}
\be\label{icbis1}
          \phi^{-1}(0) =  U\cap\partial\D, \quad \;
          \phi > 0  \;  \mbox{on $U\cap\D$ } \quad \; \mbox{and} \quad \;
          d\phi_q\not= 0  \; \mbox{for any $q\in U\cap\partial \D.$}
\ee
%\vspace*{-4.5mm}

\noindent Then,   $\partial\D$ is IC$_p$ (resp. SIC$_p$) if and only if 
for one (and hence for all) function $\phi$ satisfying (\ref{icbis1}): 

%\vspace*{-3mm}

\be\label{icbis2}
H_\phi(p)[v,v]\leq 0 \quad (\mbox{resp.} \; <0) \quad \forall v\in T_p\partial\D.
\ee

%\vspace*{-0.2mm}

(II) $\partial\D$ is {\it locally convex at} $p\in\partial \D$ (LC$_p$) if there
exists a neighborhood $U \subset  \M$ of $p$ such that

%\vspace*{-2.5mm}

\be\label{lc}
\exp_p\left( T_p\partial\D\right)\cap\left(U\cap\D\right) = \emptyset.
\ee

%\vspace*{-1.5mm}

\noindent When $\exp_p\left( T_p\partial\D\right)\cap\left(U\cap\widetilde \D\right) = \{p\}$, then 
 $\partial \D$ is {\it strictly locally   convex at}
$p$ (SLC$_p$).

\vspace{3mm}

\noindent It is easy to check 
%\vspace*{-1.5mm}
\be \label{1}
IC_p \Leftarrow LC_p, \quad \quad SIC_p \Rightarrow SLC_p.
\ee
Clearly, the converse to the first implication is not true ($\M = \R^2, 
\D= \{ (x,y): y>x^3\}, p=(0,0)$). Nevertheless, {\it if there exists a neighborhood $U$ of $p$ such that $\pD$ is IC$_q$ for all $q\in U \cap \pD$ then $\pD$ is $LC_q$ for all $q\in U\cap \D$.} Do Carmo and Warner \cite{dCW} realized that this problem is not as trivial as it sounds, and solved it (as a step for other computations) for  the constant curvature case.
Bishop \cite{Bi} solved it in general, even though differentiability $C^4$ is explicitly used in his technique. Nor does the converse to the last implication (\ref{1}) hold ($\M = \R^2, \D= \{ (x,y): y>x^4\}, p=(0,0)$).

\vspace{2mm}

\noindent The previous definitions are appliable to each point $p$ in the boundary $\pD$. The following definitions are appliable to all $\pD$:

(1)  {\it Infinitesimally   convex } (IC): 
 $\pD$ is IC$_p$  for all $p\in \D$. This is equivalent to being:

(1A) {\it locally convex} (LC):  LC$_p$, $\forall p\in \pD$ (because of Bishop's result). 

(1B) {\it variationally convex} (VC): for one (and hence for all) function $\phi: \widetilde \D \rightarrow \R$ satisfying (\ref{icbis1}) with $U=\M$,
inequality (\ref{icbis2}) holds, $\forall p \in \pD$. For the equivalence between this definition and IC, it is enough to prove the following straightforward property \cite[Cap. 3]{Bth}: for any differentiable manifold with boundary $\widetilde \D$, a function $\phi$ satisfying (\ref{icbis1}) with $U=\M$ exists.  
(VC has been widely used by using variational methods, because function $\phi$ allows direct penalization of the action functional (\ref{act}).)

The ``strict'' concepts SIC, SLC and SVC can be defined analogously, and, clearly: $SIC \Leftrightarrow SVC \Rightarrow SLC$ (but not the converse).

(2) {\it Geometrically convex} (GC):   for any $p,q\in\D$, the range of any
geodesic $\gamma:[0,1]\longrightarrow\widetilde\D$ such that
$\gamma(0) = p, \gamma(1) = q$  satisfies 
\be\label{gc}
\gamma\([0,1]\)\subset\D.
\ee
If this also holds when $p,q \in \pD$ the boundary is {\it strictly geometrically convex} (SGC). GC is a
straightforward generalization to 
Riemannian manifolds of the usual notion of convexity given
in Euclidean spaces. 
It is straightforward to check the implications $LC \Rightarrow GC \Rightarrow IC$; thus, by Bishop's result, GC is equivalent to IC. Moreover,  $ SGC \Leftrightarrow SLC$. It is worth pointing out that, in the complete case,  Germinario \cite[Theorem 2.1]{Ge}, obtained by variational methods (with a weaker assumption on differentiability than Bishop), a direct proof of $VC\Leftrightarrow GC$.

(3) {\it (Geodesically) pseudoconvex} (PC) (in the sense of \cite{BPa}, see also \cite{BEE}): for each compact set $K \subseteq \D$ there is a compact set $H \subseteq \D$ such that each geodesic segment 
$\gamma: [a,b] \rightarrow \D$ with $\gamma (a), \gamma (b) \in K$ satisfies $\gamma([a,b]) \subseteq H$. This definition is intrinsic to the (open domain of the) manifold. So, it is said that $\D$ (rather than $\pD$) is PC (thus, it is natural to assume  $\D = \M$). 

\vspace{2mm}

\noindent Easily, $ PC \Rightarrow GC$.  Nevertheless, PC is  not implied by GC; in fact, {\it a complete Riemannian manifold may be non-PC} (for example: a complete surface with infinitely many holes). Summing up: at each point $p$ implications (\ref{1}) holds and, globally

%\vspace*{-2mm}

$$PC \Rightarrow GC \Leftrightarrow LC \Leftrightarrow VC \Leftrightarrow IC.$$

%\vspace*{-1mm}

\noindent It is easy to check that, when $\widetilde \D$ is complete, there is no loss of generality assuming that so is $\M$: otherwise, a new  Riemannian metric on all $\M$ can be defined such that it is complete and coincides with the original metric on $\widetilde \D$ \cite[Cap. 3]{Bth}.
All the above equivalent conditions on convexity for $\pD$ provide different ways to prove that: 
%\vspace*{-1.5mm}
\begin{center}
{\it when $\M$ is complete, $\D$ is convex if and only if $\partial \D$ is convex.}
\end{center}
%\vspace*{-1.5mm}
Recall that the completeness assumption is essential. In fact, 
if $\M$ is incomplete and geodesically connected,  a domain $\D \subset \M$ with convex boundary may be non-geodesically connected (take as $\M$  a cylinder minus a small segment, and as $\D$ a small ball such that the segment lies inside it).
\vspace{2mm}

\noindent Now we are ready to examine the general case where $\partial \D$ is not differentiable or $\widetilde \D$ is not complete, which has been systematically studied in \cite{BGS2}. By the quoted result, it is
clear that {\it if there exists a sequence 
$\(\widetilde \D_m\), m\in\N$
of complete submanifolds with convex (differentiable) boundary such that 

%\vspace*{-2mm}

\begin{equation}
\label{seq}
\widetilde\D_m\subset \widetilde\D_{m+1} \quad \hbox{and} \quad  
\D = \bigcup_{m\in\N} \widetilde\D_m,
\end{equation}
%\vspace*{-2.5mm}

\noindent  then $\D$ is geodesically 
connected} (as a consequence, 
 the results by Gordon in  \cite{Go} are reproven and generalized in a simple way). Now, it is natural to wonder if, in this case, $\D$ must be convex. The answer is quite simple: {\it if $\widetilde \D$ is complete then $\D$ is convex; otherwise, there are counterexamples}. This result is proven by using  standard geometrical methods \cite{BGS2}.

Now, two questions arises naturally: (A) can the convexity of each $\pD_m$ be weakened? and (B)~is there any intrinsic condition for the convexity of 
$\D$ (independent of its boundary in $\M$)? To answer these questions the (intrinsic) Cauchy boundary of $\D$ and a variational approach must be taken into account. 
Let $\widetilde \D^c$ be the canonical completation of $\D$
by using Cauchy sequences, and $\partial_c \D$ the corresponding boundary points, $\widetilde \D^c = \D \cup \partial_c \D$. 
Recall that $\widetilde \D^c$ is always complete as a metric space, but the boundary points in 
$\partial_c \D$ are not necessarily differentiable and, if they are, 
the metric may be non--extendible or degenerate there. 
Note that any point of $\partial \D$ determines naturally one or more points in $\partial_c \D$, and $\widetilde \D$ is complete if and only if all the points in  $\partial_c \D$ are determined in this way by points of $\pD$. From the variational point of view, geodesics joining two fixed points $p, q \in \D$ are seen as critical points of the action functional 

%\vspace*{-2mm}

\be \label{act}
f(x) = \frac{1}{2} \int^1_0 \< \dot x(s), \dot x(s)\> ds
\ee

%\vspace*{-0.5mm}

\noindent defined on the space of differentiable curves joining $p$ and $q$ or, technically better, on its $H^1-$ Sobolev completion, i.e., the Hilbert manifold $\Omega^1(\D, p, q) = \{ x\in H^{1,2} ([0,1], \D) \mid
x(0)= p, x(1) = q\}$.

A first result \cite{BGS2} shows that a domain $\D$ as in (\ref{seq}) is convex,  even if the boundary of each $\D_m$ is not  convex (and, thus, each $\D_m$ may be non--geodesically connected), when  a suitable local estimate of the loss of convexity of $\partial\D_m$ and boundness of the sequence is ensured. 

\vspace{2mm}
\noindent {\sc Theorem 1.1.}  
Let $\(\M, \<\cdot,\cdot\>\)$ be a complete Riemannian manifold and $\D$ an open domain of $\M$. Assume that
there exists a positive differentiable function
$\phi$ on $\D$  such that

(i) $\lim_{x \rightarrow \partial \D}
\phi(x) = 0$;

(ii) each $y\in \partial \D$ admits a neighbourhood $U\subset \M$ and constants $a, b> 0$ such that 

\vspace{1mm}

\hspace{2.1cm} $a \leq \| \nabla \phi(x)\| \leq b \quad \forall x\in\D \cap U; $

\vspace{1mm}

(iii) the first and second derivatives of the normalized flow of 
$\nabla \phi$ are locally bounded close to $\partial \D$ (that is:  each $y\in \partial \D$ admits a neighbourhood $U\subset \M$ such that the induced local flow of $\nabla \phi/\| \nabla \phi \| ^2$ on $\D \cap U$ have first and second derivatives with bounded norms), and 
  
(iv) there exists a decreasing and infinitesimal  sequence 
$\{a_m\}$ such that each $y\in \partial \D$ admits a neighbourhood $U\subset \M$ and a constant $M \in \R$ satisfying:

%\vspace*{-3mm}

\be\label{1.1} 
H_\phi(x)[v,v] 
\leq M \<v,v\>\phi(x) \quad 
\forall x\in \phi^{-1}(a_m) \cap U, 
v\in T_x\phi^{-1}(a_m), m\in \N.
\ee

%\vspace*{-1.5mm}

Then $\D$ is convex.

\vspace{2mm}

\noindent It is straightforward to check that, if $\pD$ is differentiable and convex, all the conditions (i)--(iv) are automatically satisfied. 

To prove this result,  the action
functional $f$ is penalized with a term depending on a positive parameter 
$\eps$: $f_\eps (x) = f(x) + \eps \int_0^1 \phi^{-2}(x(s)) ds$. Each penalized functional is bounded from below and satisfies the condition of Palais-Smale; so, from standard variational arguments (see, for example, \cite{Ma-li}) $f_{\eps}$ admits a critical (minimum) point.
The crucial point is to prove that a critical point of $f_\eps$ in a 
sublevel of $f_\eps$ is uniformly far (with respect to $\eps$) from $\partial \D$. 
In the proof the critical points of the penalized functionals are projected (using the normalized flow of $\nabla \phi$) on the hypersurface $\phi^{-1} (a_m )$ for $m$ large enough. This makes possible to get critical points of $f$ (i.e. geodesics) not touching  $\partial \D$ by means of a limit process. 

Technical condition {\it  (iii)} and even the completeness of the ambient manifold $\M$ can be weakened if {\it (iv)} is imposed on all points and directions enough close to the boundary. So, a straightforward consequence of the technique in for Theorem 1.1 is the following result \cite{BGS2} (compare with \cite{sal1}):

\vspace{2mm}
\noindent {\sc Theorem 1.2.}
Let $\(\M, \<\cdot,\cdot\>\)$ be a Riemannian manifold, $\D \subset \M$ an open domain, 
and $\widetilde \D^c = \D \cup \partial _c\D$ its canonical Cauchy completion. Assume that
there exists a positive differentiable function
$\phi$ on $\D$  such that

(i) $\lim_{x \rightarrow \partial_c \D}
\phi(x) = 0$;

(ii) each $y\in \partial_c \D$ admits a neighbourhood $U\subset \widetilde \D^c$ and constants $a, b > 0$ such that 

\vspace{1mm}

\hspace{2.1cm} 
$a\leq \|\nabla \phi(x)\| \leq b \quad \forall x\in \D \cap U;$ 
  
\vspace{1mm}

(iii) each $y\in \partial_c \D$ admits a neighbourhood $U\subset  \widetilde \D^c$ and 
a constant $M\in \R$ such that inequality (\ref{1.1}) holds for all $x \in\D \cap U$ and for all $v \in T_x \M$.     

Then $\D$ is convex.

\vspace{2mm}
\noindent It is worth pointing out that in both previous theorems a multiplicity result can be also obtained, when the topology of the fiber is not homotopically trivial: {\it if $\D$ is not contractible (in itself), then for any $p, q\in\D$ there exists a sequence $\{x_m \}$ of geodesics in $\D$ joining them such that
$\lim_{m \rightarrow \infty} f(x_m ) = \infty.$} For the proof one uses  that the Ljusternik--Schnirelman category of $\D$ is infinite, which implies the existence of infinitely many connecting critical points of $f_\eps$ with diverging lengths.
We refer to \cite{BGS2} for a deeper discussion of these results, and for examples,  (see also  \cite{Bth}).

%\newpage

\begin{center}
\section{AFFINE CONNECTIONS}
\end{center}
\setcounter{equation}{0}

\noindent Now, consider a manifold $\M$ endowed with an affine connection $\nabla$; thus, all the properties of its geodesics will hold for all semi-Riemannian manifolds. As we are interested in geodesics, there is no loss of generality assuming that $\nabla$ is symmetric; otherwise, it is well-known that there exists another affine connection with the same geodesics and torsion-free. 

All previous notions about the convexity of the (differentiable) boundary of a domain $\D$ are directly extendible to the affine case except IC$_p$, IC, because the second fundamental form $\sigma_p$ is not canonically defined. Nevertheless, the characterization of IC$_p$ in terms of the function $\phi$ satisfying (\ref{icbis1}) does makes sense. In fact, it is not difficult to check that $(\ref{icbis2})$ still holds for a function $\phi$ satisfying (\ref{icbis1}) if and only if it holds for each such $\phi$. Thus, this will be the natural definition of IC$_p$ for affine manifolds. In principle, Bishop's result holds just in the Riemannian case and, so, the  implications which hold are (\ref{1}) and: PC or LC
$\Rightarrow$ GC $\Rightarrow$ VC $\Leftrightarrow$ IC, and 
 SLC $\Rightarrow$ SGC $\Leftarrow$ SVC $\Leftrightarrow$ SIC. 

However, these concepts (except for PC) are not so useful in the affine case, because even when $\widetilde \D$ is complete and $\pD$ is convex in the strongest sense, $\D$ may be non-geodesically connected. In fact, 
the first problem to be solved is what conditions on a manifold without boundary must be imposed to obtain geodesic connectedness. 
The following example by Bates \cite{Ba} shows that {\it a compact and complete affine 
manifold may be non-geodesically connected.}

\vspace{2mm}

\noindent {\sc Example 2.1.} Consider the moving frame $(X_1=\cos x \partial_x + \sin x \partial_y, X_2= -\sin x \partial_x + \cos x \partial_y )$ on $\R^2$, and the affine connection $\nabla$ such that $X_1, X_2$ are parallel. The geodesics of $\nabla$ are the integral curves of the (complete) vector fields $X= a X_1 + b X_2, a, b \in \R$. So, any geodesic $\gamma(s)=(x(s),y(s))$ is complete and $x(s)$ lies in an interval of length $\leq 2\pi$. Thus, a required example is the quotient torus $T^2=\R^2/4\pi\Z^2$, with the induced connection.

\vspace{2mm}

\noindent Given $p\in \M$ let Con$(p)\subseteq \M$ be the subset containing  the points of $\M$ which can be connected with $p$ by means of a geodesic. The previous example also shows that, even under strong hypotheses, Con$(p)$ may be non-closed. To study this  more in depth, consider the following concepts \cite{Lo}, \cite{BPb}.

Let $G(\M )$ be the space of the geodesics of $(\M, \nabla)$, that is, the projective bundle $P\M $ (obtained by identifying two vectors of the reduced bundle $T'\M = T\M\backslash \{ \mbox{zero section} \} $ 
if they are proportional) where two clases of non-zero vectors $[v], [w] \in P\M$ are identified if there exists a geodesic $\gamma: ]a,b[ \rightarrow \M$ such that $\gamma'(t_0) \in [v], \gamma'(t_1) \in [w]$ for some $t_0, t_1 \in ]a,b[$. Thus, each  geodesic $\gamma$ determines a unique class $[\gamma] \in G(\M )$. In the remainder of this Section, all the geodesics will be assumed {\it inextendible}. The natural quotient topology  of $G(\M )$ can be characterized as follows. Consider first  a sequence 
of geodesics $\gamma_n: ]a_n,b_n[ \rightarrow \M$ and a geodesic
 $\gamma: ]a,b[ \rightarrow \M$.  If there exists $t_0\in ]a,b[ $ contained in all but a finite number of $]a_n,b_n[$, and $\{\gamma_n'(t_0) \} \rightarrow \gamma'(t_0)$,  then lim sup$\{a_n\} \leq a < b \leq $lim inf$\{b_n\}$ and $\{\gamma'_n\}$ converges uniformly on compact subsets of $]a,b[$ to $\gamma'$ (for any distance on $T\M$ compatible with its topology) \cite[Prop. 2.1]{JMP1}. Now, a sequence of geodesics $\{\beta_n\}$ is said to {\it converge tangentially} to a geodesic $\beta$ if there exists 
$\gamma_n \in [\beta_n], \gamma \in [\beta]$ and a $t_0$ such that  
$\{\gamma_n'(t_0) \} \rightarrow \gamma'(t_0)$; this convergence holds if and only if $\{ [\beta_n]\} \rightarrow [\beta]$ in $G(\M )$, \cite{BPb}. 

\vspace{2mm}
\noindent It is worth pointing out:

(I) Tangential convergence is completely independent of completeness, even if $\M$ is compact. In fact, even if  $\gamma_n$ converges tangentially to a unique  $\gamma$, all $\gamma_n$ may be incomplete and $\gamma$ complete, or viceversa; counterexamples can be found in Lorentzian tori \cite{JMP1}.

(II) The sequence $\gamma_n$ may converge tangentially to more than one limit and, in this case, $G(\M )$ is not Hausdorff.
Remarkably, this happens when a point $p$  
of any  affine manifold 
is removed, if not all the geodesics starting at $p$ are closed. But, unfortunately, this is not by any means the only case. In fact, consider the standard flat torus $\R^2/\Z^2$; the sequence of geodesics constantly equal to $\alpha$, where $\alpha$ is induced on the torus by a geodesic in $\R^2$ with irrational slope, has infinitely many tangential limits.

\vspace{2mm}

\noindent Recall that $(\M,\nabla )$ is called {\it disprisoning} if given any geodesic $\gamma: ]a,b[\rightarrow \M$ and any compact subset $K$ of $\M$ there are sequences $\{t_n\} \rightarrow a^+$, $\{s_n\} \rightarrow b^-$ such that $\gamma(t_n)$ and $\gamma(s_n)$ do not lie in $K$. Disprisonment, pseudoconvexity and the topology of $G(\M )$ have proven to be fruitful in order to study some geometrical properties, including a Cartan-Hadamard type theorem and applications to Relativity, see \cite{Be}, \cite{BEE}, \cite{BP}. For geodesic connectedness, the following result by Beem and Parker holds
\cite{BPa}, \cite{BPb}.

\vspace{2mm}
\noindent {\sc Theorem 2.2.}
Let $(\M, \nabla)$ be  an affine manifold; the following implications are fulfilled:

$\mbox{Disprisonment and pseudoconvexity} \; \Rightarrow \; G(\M) \; \mbox{is Hausdorff} \; \Rightarrow \; \mbox{Con}(p) \; \mbox{is closed, } \; \forall p \in \M. $

\noindent Moreover, if this last property holds and there are no conjugate points, then $(\M, \nabla )$ is geodesically connected.
\vspace{2mm}

\noindent The last assertion is obvious because the absence of conjugate points implies that Con$(p)$ is open for all $p\in \M$; thus, the two implications in Theorem 2.2  
yield sufficient conditions for geodesic connectedness.

\begin{center}
\section{THE INDEFINITE SEMI-RIEMANNIAN CASE. SPACEFORMS}
\end{center}
\setcounter{equation}{0}

\noindent We refer to the standard books \cite{BEE}, \cite{O}  for definitions and general background about semi-Riemannian, and especially Lorentzian, manifolds. For  indefinite metrics, the absence of an 
associated  canonical distance   and, so, of any analog  to the Hopf-Rinow theorem, makes the problem of geodesical connectedness very subtle.  Perhaps the only non-trivial result with a clear resemblance to the Riemannian case is that of Avez \cite{Av} and  Seifert \cite{Se}: {\it in a globally hyperbolic Lorentzian manifold, any pair of causally related points {\em (i.e. which can be joined with a causal curve)} can be joined by means of a causal geodesic}. For this result is essential that, in the Lorentzian case, causal geodesics maximize locally the ``time-separation'' (or ``Lorentzian distance'') between causally related points. Global hyperbolicity introduces then a sort of compactness in the space of causal curves $\CC^{cau}_{p,q}$ joining any two fixed points $p, q$, in such a way that the lengths of curves in $\CC^{cau}_{p,q}$ are bounded, and  its supremum at every connected part of $\CC^{cau}_{p,q}$ is reached by a curve (necessarily pregeodesic) therein. 

Nevertheless, neither compactness nor completeness implies geodesic connectedness. It is interesting to study the geodesic connectedness of Lorentzian surfaces (Smith,  \cite{Sm}), in comparison with the affine case. Recall first that a plane $\PP $ endowed with a Lorentzian metric (or conformal class of metrics) is called {\it normal} if there exists a diffeomorphism of $\PP $ onto $\R^2$ which takes every null-geodesic into an axis-parallel line; this can be characterized in terms of the absence of barriers  \cite{Sm} (see also \cite{We}). Moreover, a null-complete  plane is normal if its Gaussian curvature  does not change sign off a compact subset, and the integral of its curvature is finite.

\vspace{2mm}

\noindent {\sc Theorem 3.1.} (1)  A normal Lorentzian plane is geodesically connected \cite{Sm}. (2) The universal covering of a complete Lorentzian torus with a Killing vector field $K (\not \equiv 0)$ is normal \cite{Sa-trans}.

\vspace{2mm}

\noindent Thus, a complete Lorentzian torus with such a $K$ is geodesically connected (this can be also checked more directly, \cite{Sa-trans}); the completeness in assertion (2) can be replaced by any of the following conditions (a posteriori equivalent): (i) $K$ has a definite causal sense on all the torus (timelike or null or spacelike), or (ii) the metric is (globally) conformally flat.  
Now, consider the following Lorentzian torus \cite[Sect. 5]{Sm} (see also \cite[Sect. 3]{Be}, 
\cite[Sect. 5]{Sa-trans}). 

\vspace{2mm}
\noindent {\sc Example 3.2.} Take on $\R^2$ the moving frame $X_1, X_2$ as in Example 2.1, and consider the Lorentzian metric $g$ such that $X_1$ and $X_2$ are null and $g(X_1,X_2)=-1$. The velocities of any timelike or spacelike curve $(x(s),y(s))$ must remain in any of the four cones continuously determined  by $X_1, X_2$; so, the length of the projection of $x(s)$ is  again bounded. Thus, a geodesically disconnected Lorentzian torus is induced.

\vspace{2mm}

\noindent It is worth pointing out that in this example $K=\partial_y$ is a Killing vector field (in Example 2.1 $K$ is an affine vector field) inducible on the torus. So, this Lorentzian torus is geodesically incomplete (the causal character of $K$ changes); we do not know any example of complete and geodesically disconnected torus. 

\vspace{2mm} 
\noindent The fact that a complete semi-Riemannian manifold may be geodesically disconnected can be stressed studying spaceforms. We say that a semi-Riemannian $n$-dimensional manifold $\M$ of index $\nu$  is a spaceform if it is complete with constant curvature $C$. In this case, $\M$ is covered by the corresponding model (1-connected) space $M(n,\nu,C)$, that is, $\M=M(n,\nu,C)/\Gamma$, where $\Gamma $ is the fundamental group of $\M$. The case $C=0$ is trivial (the model space $M(n,\nu,0)\equiv \R^n_{\nu}$ is geodesically connected) and, up to a homothety, we can assume $C=1$ (the homothetic factor may be positive as well as negative; recall that  the  Levi-Civita connection remains unchanged). If $n\geq 3$, the model space is then the pseudosphere    $\S^n_{\nu}$ (spacelike vectors of norm 1 in $\R^{n+1}_{\nu}$); the Lorentzian pseudosphere $\S^n_{1}$ is also called {\it de Sitter} spacetime, and it is globally hyperbolic. Recall that an affine manifold is called {\it starshaped} from a point $p$ if $\exp_p$ is onto. It is not difficult to check \cite[Propos. 5.38]{O} that {\it no indefinite pseudosphere $\S^n_\nu$, 
$0<\nu < n$ is geodesically connected.} The following result by Calabi and Markus \cite{CM}, in particular, solves completely the geodesic connectedness of Lorentzian spaceforms of positive curvature with $n\geq 3$.

\vspace{2mm}
\noindent {\sc Theorem 3.3.} For  $n\geq 2$:

(1) Two points $p, q \in \S^n_1$ are connectable by a geodesic if and only if $\<p,q\>_1 >-1$,  where $\<\cdot,\cdot \>_1$ is the usual Lorentzian inner product of $\L^{n+1} \equiv \R^{n+1}_1$.

(2) Every spaceform $\M= S^n_{1} /\Gamma, \M\neq S^n_{1}$ is starshaped 
from some point $p\in \M$.

(3) A spaceform $\M= S^n_{1}/\Gamma$ is geodesically connected if and only if it is not time-orientable.

\vspace{2mm}

\noindent  The proof of (1) follows by a direct computation of the geodesics. For the remainder, it is essential that, whenever $2\nu \leq n$, the  group of isommetries $\Gamma$ is finite. Then, up to conjugacy, $\Gamma \subset O(1)\times O(n) \subset O_1(n+1)$, and the proof follows by studying the barycenter of the orbits, which must lie in the timelike axis of $\R_1^{n+1}$. 

For arbitrary index $\nu$  (including the case $\nu = n-1$, which is 
equivalent to the Lorentzian case of constant negative curvature) the only general results we know are extensions of Theorem 3.3, with $n\geq 3$, and  assuming as an additional hypothesis (when $2\nu >n$) that the fundamental group $\Gamma$ is  finite. For these extensions, the non time-orientability must be replaced by the inexistence of a {\it proper time-axis} \cite[Theor. 11.2.3]{Wo}. We recall that a {\it time-axis} $T$ is a one-dimensional $\Gamma$-invariant negative-definite linear subspace of $\R^{n+1}_{\nu}$. $T$ is {\it proper} if $\Gamma$ acts trivially on $T$; that is, if $T$ is not proper then $\Gamma$ also acts as a multiplication by -1.

\vspace*{0.5mm}

\begin{center}
\section{VARIATIONAL METHODS.}
\end{center}
\setcounter{equation}{0}

\noindent The systematic application of variational methods on infinite-dimensional manifolds to Lorentzian Geometry started with a seminal paper by Benci and Fortunato \cite{BF} (see also \cite{BFa}), who studied geodesic connectedness of standard stationary spacetimes. Since then, the geodesic connectedness of stationary as well as splitting spacetimes, with or without boundary, has been widely studied by variational methods. We will review briefly these results, explaining mainly the stationary case, and giving some references and comments for the splitting case  (see also the next Section). We refer  to \cite{Ma-li} for  general background on  variational methods and applications to this and other problems in Lorentzian Geometry, and to \cite{Sa-gre} for properties of Killing fields which  will be borne in mind.

\vspace{2mm}
\noindent For stationary manifolds we will follow the approach in \cite{GP}, because it makes more intrinsic assumptions and  recovers   the previous results. Recall that a Lorentzian manifold $(\M,\<\cdot,\cdot\>)$ is called {\it stationary} %(spacetime) 
if it admits a globally defined timelike Killing vector field $K$. Fixing two points $p, q \in \M$,  geodesics joining them are, still in the  semi-Riemannian case, the critical points of the action functional $f$  on $\Omega^1(\M, p, q)$ given in (\ref{act}). But now this functional is bounded neither from above nor from below, and it may not satisfy Palais-Smale condition. Nevertheless, this problem can be solved by taking into account that any critical point  of $f$ (or geodesic) $z$ satisfies $\<K,z'\> \equiv C_{z}$, where $C_{z}$ is a constant. This suggests that variations in the $K$--direction are irrelevant, and that, under the orthogonal splitting of the tangent bundle $TM= \mbox{Span}(K) \oplus  \mbox{Span}(K)^{\perp}$, only projections on the spatial part $\mbox{Span}(K)^{\perp}$ are important  (similar ideas work for geodesic completeness  \cite{RS-pams}, 
\cite{RS-gd}). %, \cite{Sa-gre}

More precisely, let $\CC_{p,q}=\{ z\in C^1([0,1], \M) \mid z(0)=p, z(1)=q, \<K,z'\> \equiv C_{z}\}$ and $\NN_{p,q}$ its Sobolev $H^{1,2}$-completation; for technical computations, the Riemannian metric $g_R$ obtained by reversing $\<\cdot, \cdot\>$ on Span$(K)$ can be used.
The action functional $f$ can be written as the sum of two functionals $f_1, f_2$, 

%\vspace*{-2mm}

$$f_1(z) =  \frac{1}{2} \int^1_0 \left( \< \dot z(s), \dot z(s)\> 
- \frac{\< \dot x(s), K\>^2}{\<K,K\>} \right)ds, \quad 
f_2(z) = \frac{1}{2} \int^1_0 \frac{\< \dot z(s), K\>^2}{\<K,K\>} ds. $$

%\vspace*{-1mm}

\noindent It is not difficult to check that $f'_1$ vanishes on all tangent vectors $\zeta$ to 
$\Omega^1(\M, p, q)$ which are pointwise parallel to $K$, that is, $f_1$ is the ``spatial'' part with respect to $K$ of the functional $f$. Moreover, $\NN_{p,q}$ can be characterized as the set of $z\in\Omega^1(\M, p, q)$ such that $f'(z)[\zeta] = f'_2(z)[\zeta]=0$, for any such $\zeta$ parallel to $K$; by the implicit function theorem, $\NN_{p,q}$ is a submanifold of $\Omega^1(\M, p, q)$. 
Finally, one can check that 
critical points of $f$ on $\Omega^1(\M, p, q)$ coincide with the critical points of the restriction $J$ of $f$ to  $\NN_{p,q}$ (and thus to  $\CC_{p,q}$). 

Summing up, {\em if  some conditions are imposed on $\CC_{p,q}$  so that the restriction $J$ of $f$ to $\NN_{p,q}$ must reach a critical point, then a geodesic connecting $p$ and $q$ is obtained.} 
And the most natural variational such conditions are: 

%\vspace*{-2mm}

\begin{quote}
{\it $\CC_{p,q}$ is (i) non-empty and (ii) $c-$precompact, for some $c>{\rm Inf}_{\CC_{p,q}}f(z)$} 
\end{quote}

%\vspace*{-2mm}

\noindent ($\CC_{p,q}$ is said {\it $c-$precompact} if every sequence $\{z_n\} \subset \CC_{p,q}$ with $f(z_n)\leq c$ has a uniformly convergent subsequence  in $\M$). Thus, if (i) and (ii) hold for all $p,q \in \M$, then $(\M, \<\cdot, \cdot \>)$ is geodesically connected (one can check that 
$J$ is bounded from below, its sublevels $J^{c'}$ are complete metric subspaces of $\NN_{p,q}$,  for all $c'\leq c$, and 
Palais-Smale condition is fulfilled). 

Now, we can wonder when (i) and (ii) hold. 
If the restriction of $f$ to $\CC_{p,q}$ is {\it pseudocoercive} (that is, 
$c-$precompact for all $c\geq {\rm Inf}_{\CC_{p,q}}f(z)$) for any $p, q \in \M$
then $(\M,\<\cdot , \cdot \>)$ is globally hyperbolic, but the converse is not true. Condition (i) holds if either $p$ and $q$ are causally related or  $K$ is complete (this happens, for example, when the auxiliary Riemannian metric $g_R$ is complete, and $g_R(K,K) (= -\<K,K\>)$ is bounded); clearly, the converse does not hold.  Assume that $(\M,\<\cdot , \cdot \>)$ is a {\it standard stationary spacetime}, that is, 
$\M$ is a product manifold 
$\M= \R \times  \M_0 $ ($\M_0$ any manifold) and $\< \cdot , \cdot \>$ can 
be written, with natural identifications, as:

%\vspace*{-1.5mm}

\begin{equation} \label{st-st}
\<\cdot , \cdot \>
= -\beta dt^{2} +  2 \omega \otimes dt  + g_0, 
\label{st}
\end{equation}

%%\vspace*{-1mm}

\noindent where 
$dt^{2}$ is the usual metric on $\R$, 
and $g_0,\ \beta,\  
\omega$ are, resp., a Riemannian metric, 
a positive function and a 1-form, all on $\M_0$ (locally,  stationary spacetimes look like  standard stationary ones). 

\vspace{2mm}

\noindent {\sc Theorem 4.1.}
A standard stationary spacetime is geodesically connected, if: (a) $g_0$ is complete, (b) $0<$Inf$(\beta) \leq$ Sup$(\beta) < \infty$, and (c) the $g_0$-norm of $\omega(x)$ has a sublinear growth in $\M_0$. 

\vspace{2mm}

\noindent (For (c), we mean that the norm of $\omega(x)$ has an upper bound $A\cdot d_0(x,p_0)^\alpha + B$, for some $A, B \in \R, \alpha \in [0,1[, p_0 \in \M_0$, where $d_0$ is the $g_0$-distance on $\M_0$). In fact, under these three conditions, assumptions (i) and (ii) are always satisfied;
we refer to \cite[Prop. A.3]{GP} for a more intrinsic way to express these conditions, in terms of stationary manifolds admitting a differentiable time function (which are standard stationary {\it a posteriori}). In the standard static case ($\omega \equiv 0$) condition $0<$Inf$(\beta)$ can be dropped (see \cite{Bth}, \cite{Ba-jde}); however, we remark that the imposed inequalities always imply global hyperbolicity \cite[Cor. 3.4, 3.5]{Sa-ba}.

\vspace{2mm}
\noindent It is also worth pointing out:

(I) This technique provides also consequences for the existence of infinitely many connecting geodesics or timelike geodesics when $\M$ is not contractible. In fact, 

{\it  two points $p, q$ of a stationary manifold can be joined by a sequence of spacelike geodesics with divergent lengths if $K$ is complete, $\CC_{p,q}$ is pseudocoercive and $\M$ is non-contractible}; 

\noindent (the essential step for the proof is that $\Omega_{p,q}(M)$ is homotopically equivalent to $\NN_{p,q}$ and, thus, the Ljusternik-Schnirelman category of $\NN_{p,q}$ is infinite).

For timelike geodesics, recall that, under our type of assumptions, Avez-Seifert's  technique is appliable and chronologically related points can be joined by timelike geodesics. Let $p,q \in \M$, and $\gamma_q(t)$ be an integral curve of $K$ starting at $q$; when $K$ is complete then $p$ belongs to the chronological past of $\gamma_q(t)$ for $t$ big enough (a direct proof is not difficult, see also \cite[Sect. 4]{Sa-ba}).  Then, under precompactness,  $\NN_{p, \gamma_q(t)}$ contains at least a timelike geodesic for $t$ big enough and, when the topology is not homotopically trivial, the following result on multiplicity holds (see \cite[Theorem 1.4]{GP}):

{\it If $\M$ is non-contractible, $K$ is complete and there exist $c_0 <0, t_0>0$ such that $\NN_{p, \gamma_q(t)}$ is $c_0$ precompact for all $t>t_0$, then the number of timelike geodesics joining $p$ and $\gamma_q (t)$ goes to $\infty$ when $t\rightarrow \infty$.}

(II) Analogous techniques should work if a semi-Riemannian manifold $(\M,g)$ of index $s$ admits $s$  Killing vector fields $K_1, \dots ,K_s$ independent at each $p\in \M$ such that $g$ restricted to $\KK =$ Span$\{K_1, \dots , K_s \}$ is negative definite. We have now the natural splitting $TM = \KK \oplus \KK^\perp$, and, for each geodesic $z$,  the projection of $\dot z$ on $\KK$ can be recovered from the constants $C_{z,i} \equiv g(\dot z, K_i)$  \cite{GPS} (the analogous problem for geodesic completeness was  solved  in \cite{RS-gd}). 

Moreover, even in the Lorentzian case, one can consider the case when there exist two pointwise independent Killing vector fields $K_1, K_2$ such that  
$\{ K_1(p), K_2(p)\}$
spans a Lorentzian plane at each $p\in \M$ but  neither $K_1$ nor $ K_2$ are timelike on (all) $\M$.  Remarkably, this happens in G\"odel type spacetimes; for the modifications of the technique in this case, see \cite{CS}, \cite{Ca}. 

(III) Let us discuss the case with boundary briefly (see also, 
\cite{Bth}, \cite{BS}, \cite{Ge2}). Consider first a standard stationary spacetime $\M = \R \times \M_0 $ satisfying the assumptions of Theorem 4.1  (thus, geodesically connected), and let $\D_0 \subset \M_0$ be an open domain with differentiable boundary $\pD_0$. We have seen that, in general, a domain $\D$ of a complete or geodesically connected semi-Riemannian manifold does not inherit good properties for geodesic connectedness. Nevertheless, if $\D = \R \times \D_0 \subset \M$ is variationally convex (VC), then the corresponding function $\Phi$ can be chosen  independent of $t$, and the functionals $f$,  $J$, can be penalized in a similar way to the Riemannian case. So, penalized functionals $J_\epsilon$ do satisfy Palais-Smale condition, and one obtains geodesic connectedness.  Moreover, one obtains again that $\pD$ is VC if and only if it is GC \cite{BGS3}. Of course, if we are interested just in $\D$, it is not exactly  relevant for $\M$ to fulfill assumptions of Theorem 4.1: 
one needs just
the convexity of $\pD$ and the possibility to extend the metric on $\D$ to all $\M$, in such a way that the assumptions of Theorem 4.1 are fulfilled; this possibility can be expressed as more intrinsic conditions  on $\D$. Furthermore, the assumption that the stationary manifold is standard can also be dropped. In fact, in order to penalize the functionals $f, f_1, f_2$ as in the Riemannian case, one  needs only that the boundary $\pD$ can be determined by a function $\phi$ as in the definition  of VC, which is also invariant by the flow of $K$.

When the boundary is not smooth, the problem has been studied in the standard {\it static} case ($\omega \equiv 0$). A result in the spirit of Theorem 1.2 can be obtained, proving the geodesic connectedness of spacetimes such as outer Reissner-Nordstr\"om's and Schwarzschild's \cite{BFG-st} (see also \cite[Cap. 6]{Bth}). On the other hand, for the standard static case, a different variational principle in \cite{Sa-nonlin} solves completely the problem of connecting a point and a integral curve of $\partial_t$, for 
arbitrary $\M_0$ (or $\D_0$), \cite{BGS3}.

\vspace{2mm}

\noindent We will mean by a {\it splitting spacetime} a product manifold $\M = \R \times \M_0$ endowed with a Lorentzian metric as (\ref{st-st}) but allowing  
 $g_0,\ \beta,$ and   
$\omega$ to depend  differentiably on the time variable $t$. In this case,  bounds on the derivatives with respect to $t$ of these elements  must be also imposed in order to obtain geodesic connectedness. As a typical result we have \cite{GM-ihp}:

\vspace{2mm}

\noindent {\sc Theorem 4.2.}
A splitting spacetime $(\M, \< \cdot, \cdot \> =  -\beta (t,x) dt^{2} +  2 \omega (t,x)\otimes dt  + g_t(x))$ is geodesically connected, if: 

(1) $g_0$ is complete and there exists $\lambda>0$ such that $g_t > \lambda g_0$ for all $t$.

(2) $0<$Inf$(\beta) $ and $\beta(x,0), \omega(x,0)$ are bounded. 

(3)  $g_t/\beta(t,x)$ and $\omega/\beta(t,x)$ are bounded by a function on $\M$ type: $b_0(x) + b_1(x) |t|^{\mu},$ ($\mu \in [0,1[$, and $\mu \in [0,2[$, resp.) 

(4) Consider  the natural derivatives $\partial_t \alpha , \partial_t \beta, \partial_t \delta$ of $\alpha, \beta, \delta$ with respect to $t$, resp. Then the $g_t$--norms of 
$\partial_t \alpha/\alpha , \partial_t \beta /\beta, \partial_t \delta$ are bounded at each hypersurface with constant $t$, and its supremum when 
$t \rightarrow \pm \infty$ goes to 0. 

\vspace{2mm}

\noindent Moreover, domains of type $\D = ]a,b[ \times \M_0 $ ({\it strips}) are shown to inherit geodesic connectedness, provided that $\pD$ is VC (see \cite{BFM} for orthogonal splittings $\omega \equiv 0$, 
\cite{GM-ihp} for non-orthogonal splittings, and \cite{Ma2} for strips; see also \cite{Ma-li}).  

Nevertheless, now the $t$-dependence does not allow a reduction to an equivalent Riemannian problem, as in the stationary case. Instead,  Rabinowitz's saddle point theorem \cite{Ra} is used, but  two technical complications must be circumvented: (A) $f$ does not satisfy a Palais-Smale condition, which is solved by approximating  by a family of functionals $f_\eta, \eta \geq 0, \ f_0 = f $, and making some a priori estimates of the critical points of $f_\eta, \eta >0$ in order to ensure a good behavior under the limit $\eta \rightarrow 0$, and (B) in Rabinowitz's theorem, the independent directions where the functional goes to $-\infty$ are finite; so, a Galerkin finite-dimensional approximation is carried out. For the existence of infinitely many connecting geodesics when $\M$ is not contractible, the {\it relative category}, a  topological invariant somewhat subtler then the Ljusternik-Schnirelman category, is used.

Finally, it is worth pointing out that the problem of geodesic connectedness  is naturally generalized to others like:
(1) the  connectedness of two submanifolds by normal geodesics, studied for splitting manifolds in \cite{CMS}, or (2) the 
 connection by trajectories of some more general Lagrangian systems, studied for stationary manifolds and potential vector fields independent of time in \cite{Ba-jde}.

\begin{center}
\section{MULTIWARPED SPACETIMES. A TOPOLOGICAL METHOD}
\end{center}
\setcounter{equation}{0}

\noindent A {\it multiwarped spacetime} is a  product manifold $I\times F_1 \times \cdots \times F_m,$ $I= ]a,b[ \subseteq \R{}$   endowed with a  metric $g = -dt^2 + \sum_{i=1}^m f_i^2(t) g_i, t\in I$, where $f_1, \dots f_m$ are positive functions on $I$, and each $g_i$ is a Riemannian metric on the manifold $F_i$. These spacetimes include classical examples of spacetimes: when  $m=1$ they are the {\it Generalized Robertson-Walker} (GRW) spacetimes, standard models of inflationary spacetimes \cite{Sa98}; when $m=2$, the intermediate zone of Reissner-Norsdstr\"om spacetime and the interior of Schwarzschild spacetime appear as particular cases \cite{Sa97}; moreover, multiwarped spacetimes may also represent relativistic spacetimes together with internal spaces attached at each point (see \cite{Sc}). 

The geodesic connectedness of this type of spacetimes with $m=2$ have been studied by using variational methods in manifolds  without and with   boundary  \cite{Gi}, \cite{GM-mm}. Nevertheless, more accurate results are obtained  by using a topological method introduced in \cite{FS-trans}; in fact, the following result is proven:

\vspace{2mm}
\noindent {\sc Theorem 5.1.} (1) In a multiwarped spacetime with convex fibers $(F_{1},g_1),\dots, (F_{m},g_m)$  each two causally related points can be joined by a causal geodesic. 

\vspace{1mm}

(2) The multiwarped spacetime is geodesically connected if  the fibers are and:
 
%\vspace*{-3mm}
\begin{equation}\label{e28} 
\begin{array}{ll} 
\int_{c}^{b}f_{i}^{-2}(f_{1}^{-2}+\cdots +f_{m}^{-2})^{-1/2}=\infty  & \int_{a}^{c}f_{i}^{-2}(f_{1}^{-2}+\cdots +f_{m}^{-2})^{-1/2}=\infty  
\end{array} 
\end{equation}
%\vspace*{-4.5mm}

\noindent for all $i$ and for some $c\in (a,b)$. 
Moreover, if one of the fibers $F_j$ is not contractible then: 

(a) each two points can be joined by infinitely many geodesics, and 

(b) for any $z\in I\times F_1 \times \cdots \times F_m,$ and $x \in 
 F_1 \times \cdots \times F_m,$
the number of timelike geodesics joining  $z$ and $(t,x)$ goes to $\infty$
when $t$ goes to an endpoint of the interval $I$.

\vspace{3mm}

\noindent It is worth pointing out:

(I) Equality (\ref{e28}) is equivalent to the following condition: any point $z_0$ of the spacetime can be joined with any line $L[x]$ by means of both, a future directed and a past directed causal curve. 
Nevertheless, Theorem 5.1 does not cover all the possibilities of the technique,  and more general versions of this theorem can be given.
These general versions give very accurate results; in fact,  a necessary and sufficient condition for geodesic connectedness when $m=1$ can be given with a reasonably long distinction of cases \cite{FS-jgp} (multiplicity, existence of timelike geodesics, conjugate points, etc. are also completely characterized in this reference; see also \cite{Sa98}). 
In particular, geodesic connectedness of Reissner-Nordstr\"{o}m Intermediate spacetime is reproven (previous proofs and results in this direction were obtained in \cite{GM-mm}, \cite{Sa97}). The accuracy of the technique is shown by proving the geodesic connectedness of Schwarzschild inner spacetime; in fact,  a good behaviour of the warping functions yields geodesic connectedness, but  the warping functions of Schwarzschild inner spacetime do not have such good behaviour. Nevertheless, in this case, this problem can be skipped because one of the fibers of Schwarzschild's is a sphere (and so, each pair of its points can be joined by geodesics of arbitrarily large length); if  this fiber is replaced  by a plane, the resulting spacetime is {\it not}  geodesically connected.

(II) If the fibers are assumed to have boundary, the problem is reduced to considering  when this boundary implies geodesic connectedness or convexity (Section 1). The general technique also works when a strip $]a',b'[ \subset I$ is considered; so, the problem with boundary is also solved, for boundaries which preserve the multiwarped structure.

\vspace{2mm}

\noindent For the proof of Theorem 5.1, the Avez-Seifert type result (1) relies on a partial integration of the geodesics. For (2), the idea is the following, under an assumption somewhat stronger than (\ref{e28}) (under (\ref{e28}) some technicalities must be also taken into account).
As  there are points in each $L[x]$ both, future and past related with $z$,
 these points can be joined with causal geodesics. Thus, we have just to connect $z=(z_0,z_1,\dots z_m)$  and $(t,x), x=(x_1,\dots x_m)$ for $t$ in  a compact interval. For this: (A) fixing the geodesics in the fibers joining each $z_i$ and $x_i$,  
a  continuous map $\bar \mu(c,K)$,  $\bar \mu: ]0,1[^{m-1} \times [K^-,K^+] \longrightarrow \R^{m-1}$, is constructed in such a way that each zero of $\bar \mu$ represents the initial condition of a connecting geodesic, (B)  for 
$K= K^-$ and $K=K^+$, the result above on causal geodesics ensures the existence of   at least one zero of (a continuous  extension of) $\bar \mu$  for some $c\in [0,1]^{m-1}$; then, the problem is solved if a continuous set of zeros ${\cal Z}$ containing these two zeros is found, (C) for each $K \in ]K^-,K^+[$  the behavior of $\bar \mu(c,K)$ at the boundary of $[0,1]^{m-1}$, 
allows to find ${\cal Z}$, by means of arguments  on continuity of solutions of equations depending on a parameter, based in Brower's 
topological degree. When $F_j$ is not contractible, the result follows by applying the technique for each geodesic joining $z_j$ and $x_j$.


\newpage




\begin{thebibliography}{99}



{\ixpt 


\bibitem  {Av} AVEZ A., 
Formule the Gauss-Bonnet-Chern en m\'etrique de signature quelconque, 
{\it C.R. Acad. Sci.} {\bf 255,} 2049-2051 (1962).

%\vspace{-1.5mm}

\bibitem{Bth} BARTOLO R., {\it Curvas cr\'{\i}ticas en variedades riemannianas y lorentzianas con borde}, Ph.D thesis, Univ. Granada, 2000.

%\vspace{-1.5mm}

\bibitem{Ba-jde} BARTOLO R.,
 Trajectories connecting two events of a Lorentzian manifold
in the presence of a vector field, 
{\it J. Diff. Equat.}
{\bf 153}, 82-95 (1999)

%\vspace{-1.5mm}

\bibitem{BGS1} BARTOLO R., GERMINARIO A. \& S\'ANCHEZ M.,
 Periodic trajectories with fixed energy on Riemannian and Lorentzian
manifolds with boundary,
{\it Ann. Mat. Pura Appl. (IV)} {\bf CLXXVII}, 241-262 (1999)

%\vspace{-1.5mm}

\bibitem{BGS2} BARTOLO R., GERMINARIO A. \& S\'ANCHEZ M., Convexity of domains of Riemannian manifolds, preprint.

%\vspace{-1.5mm}

\bibitem{BGS3} BARTOLO R., GERMINARIO A. \& S\'ANCHEZ M.,
 A note  on the boundary of  a static spacetime,  
preprint.

%\vspace{-1.5mm}

\bibitem{BS} BARTOLO R. \& S\'ANCHEZ M., these proceedings.


\bibitem{Ba} BATES L., 
You can't get there from here, 
{\it Diff. Geom. Appl.}, {\bf 8} No. 3, 273-274, (1998).

%\vspace{-1.5mm}

\bibitem{Be} BEEM J.K. 
Disprisoning and pseudoconvex manifolds, 
{\it Proc. Symp. Pur. Math.} {\bf 54} 19--26 (1993) Part 2.

%\vspace{-1.5mm}

\bibitem{BEE}  BEEM  J.K.,  EHRLICH P.E. \& EASLEY K., {\it Global  Lorentzian  Geometry,}   %Pure  Appl. Math., 
Marcel Dekker Inc., N.Y. (1981).

%\vspace{-1.5mm}

\bibitem{BP}  BEEM  J.K. \& PARKER P.E. 
Pseudoconvexity and general relativity, 
{\it J. Geom. Phys.}  {\bf 4}, 71--80 (1987).

%\vspace{-1.5mm}

\bibitem{BPa}  BEEM  J.K. \& PARKER P.E. 
Pseudoconvexity and geodesic connectedness, 
{\it Ann. Mat. Pura Appl.}  {\bf 155}, 137--142 (1989).

%\vspace{-1.5mm}

\bibitem{BPb}  BEEM  J.K. \& PARKER P.E. 
The space of geodesics, 
{\it Geom. Dedicata} {\bf 38}, 87--99 (1991).

%\vspace{-1.5mm}

\bibitem  {BF} BENCI V. \& FORTUNATO D., 
Existence  of  geodesics for the Lorentz metric of a stationary gravitational field, 
{\it Ann. Inst. Henri Poincar\'e'} {\bf 7,}  27-35 (1990).

\bibitem  {BFa}  BENCI V. \& FORTUNATO D., Periodic trajectories for the Lorentz  metric of  a  static  gravitational  field,  {\it Proc.  on Variational   Methods}  (Edited by H. BERESTICKY, J.M. CORON and I. EKELAND) 413-429, Paris (1988).

%\vspace{-1.5mm}

\bibitem  {BFG-st} BENCI V., FORTUNATO D. \& GIANNONI F., 
{\em On the existence of multiple geodesics in static spacetimes}, 
{\it Ann. Inst. H. Poincar\'e}, {\bf 8}, (1991) 79-102.

%\vspace{-1.5mm}

\bibitem  {BFM} BENCI V., FORTUNATO D. \& MASIELLO A., On the geodesic connectedness of Lorentzian manifolds,
{\it Math. Z.} {\bf 217}, 74-94 (1994).
%\vspace{-1.5mm}

\bibitem{Bi}  
BISHOP R.L. 
Infinitesimal convexity implies local convexity, 
{\it Indiana Math. J.} {\bf 24} No. 2,   169-172, (1974).

%\vspace{-1.5mm}

\bibitem{CM}  CALABI E. \&  MARKUS L.,  
Relativistic space forms, 
{\rm Ann. Math.}, {\bf 75}, 63-76, (1962).

%\vspace{-1.5mm}


\bibitem {CMS} CANDELA A.M., MASIELLO A. \& SALVATORE A., 
Existence and multiplicity of normal geodesics in Lorentzian manifolds, 
{\it J. Geom. Anal.}, to appear.

%\vspace{-1.5mm}

\bibitem {CS} CANDELA A.M. \& S\'ANCHEZ M., 
Geodesic connectedness in G\"odel-type spacetimes, 
{\it Diff. Geom. Appl.}, to appear.

%\vspace{-1.5mm}
\bibitem {Ca} CANDELA A.M. \& S\'ANCHEZ M., these proceedings.


\bibitem{dCW}
DO CARMO M.P. \& WARNER F.W., 
 Rigidity and convexity of hypersurfaces in spheres, {\it J. Differ. Geom.} 
{\bf 4}, 133--144 (1970)

%\vspace{-1.5mm}

\bibitem{FS-trans}
FLORES J. \& S\'ANCHEZ M.,
 Geodesic connectedness of multiwarped spacetimes, preprint.

%\vspace{-1.5mm}

\bibitem{FS-jgp}
FLORES J. \& S\'ANCHEZ M., 
 Geodesic connectedness and conjugate points in GRW spacetimes, {\it J. Geom. Phys.}, to appear.

%\vspace{-1.5mm}

\bibitem{Ge}  
GERMINARIO A.,  
Homoclinics on Riemannian manifolds with convex boundary, 
{\it Dynam. Syst. Appl.}, {\bf 4}, 549-566, (1995).

%\vspace{-1.5mm}
\bibitem{Ge2} GERMINARIO A., these proceedings.



\bibitem{Gi}
GIANNONI F.,
	Geodesics on nonstatic %Lorentz
 manifolds of Reissner--Nordstr\"om 	type, 
   {\it Math. Ann.} {\bf 291}, 383-401 (1991).

%\vspace{-1.5mm}

\bibitem{GM-mm}   
GIANNONI F. \& MASIELLO A.,
Geodesics on Lorentzian manifolds with quasi--convex boundary, 
{\it Manus. Math.} {\bf 78}, 81-96 (1993)

%\vspace{-1.5mm}

\bibitem{GM-ihp}   
GIANNONI F. \& MASIELLO A.,
 Geodesics on product Lorentzian manifolds, 
{\it Ann. Inst. H. Poincar\'e, %Anal. Non Lin\'eaire
} {\bf 12}, 27-60 (1995)

%\vspace{-1.5mm}

\bibitem{GP}  
GIANNONI F. \& PICCIONE P.,  
An intrinsic approach to the geodesical connectedness of stationary Lorentzian manifolds,
 {\it Comm. Anal. Geom.}, {\bf 7}, 157--197, (1999).

%\vspace{-1.5mm}

\bibitem{GPS}  
GIANNONI F., PICCIONE P. \& SAMPALMIERI R.,   
On the geodesical connectedness for a class of semi-Riemannian manifolds,
{\it J. Math. Anal. Appl.}, to appear.

%\vspace{-1.5mm}

\bibitem{Go}   
 GORDON W.B.,  
The existence of geodesics joining two given points 
  {\it J. Diff. Geom.} {\bf 9}, 443-450, (1974).

%\vspace{-1.5mm}

\bibitem{Lo} LOW R.J.,  
The geometry of the space of null geodesics
{\it J. Math. Phys.} {\bf 30}, 809--811, (1989).

%\vspace{-1.5mm}

\bibitem  {Ma-li} MASIELLO A., {\it Variational methods in Lorentzian Geometry,} %Pitman Res. Notes  Math. Ser. {\bf 309} 
Longman Sc. Tech., Harlow, Essex (1994).

%\vspace{-1.5mm}

\bibitem  {Ma2} MASIELLO A., Convex regions of Lorentzian manifolds, {\it Ann. Mat. Pura Appl. (IV)} {\bf CLXVII}, 299-322 (1994).

%\vspace{-1.5mm}

\bibitem{Sc} MIGNEMI S. \& SCHMIDT H.-J.: Classification of multidimensional inflationary models, {\it J. Math. Phys.} {\bf 39} No. 2 (1998) 998-1010.



\bibitem {O} O'NEILL B., {\it Semi-Riemannian Geometry, 
%with  applications  to  Relativity
} Ser.  Pur.  Appl. Math. {\bf 103} Academic Press, N.Y. (1983).

%\vspace{-1.5mm}

\bibitem{Ra}
RABINOWITZ P.H.,
{\em Min--Max methods in critical point theory and the min-max principle}, 
CBMS Reg. Conf. Soc. in Math. n.65: AMS, 1984.

%\vspace{-1.5mm}

\bibitem{JMP1} ROMERO A. \& S\'ANCHEZ M., 
On the completeness of 
geodesics obtained as a limit 
{\it J. Math. Phys.} {\bf 34}, 3768--3774, (1993).

%\vspace{-1.5mm}

\bibitem  {RS-gd} ROMERO A. \&  S\'ANCHEZ M., 
On the completeness of certain families  of semi-Riemannian manifolds, 
{\it Geom. Dedicata} {\bf 53,} 103-117 (1994).

%\vspace{-1.5mm}

\bibitem  {RS-pams} ROMERO A. \&  S\'ANCHEZ M.,  
Completeness of compact  Lorentz  manifolds admitting a timelike conformal Killing vector field, 
{\it Proc.  Amer.  Math.  Soc.} {\bf 123,}  2831-2833 (1995).

%\vspace{-1.5mm}

\bibitem{sal1}  SALVATORE A.,  
A two points boundary value problem on non complete Riemannian manifolds, 
	{\it Variational Methods in Nonlinear Analysis}, A. Ambrosetti, K.C. Chang Eds., Gordon and Breach Pub.,    149-160, (1995). 

%\vspace{-1.5mm}

\bibitem{Sa-nonlin}
S\'ANCHEZ M.,
	Geodesics in static spacetimes and $t$--periodic trajectories, 
     {\it Nonlin. Anal.} {\bf 35}, 677-686 (1999).

%\vspace{-1.5mm}

\bibitem{Sa98}
S\'ANCHEZ M.,
	On the geometry of 
%GRW 
Generalized Robertson Walker spacetimes,
   {\it Gen. Relat. Gravit.} {\bf 30}, 915-932 (1998)

%\vspace{-1.5mm}

\bibitem{Sa97}
S\'ANCHEZ M.,
	Geodesic connectedness in Generalized Reissner--Nordstr\"om type 	Lorentz manifolds, 
  {\it Gen. Relat. Gravit.} {\bf 29}, 1023-1037 (1997)

%\vspace{-1.5mm}

\bibitem  {Sa-trans} S\'ANCHEZ M., 
Structure of Lorentzian tori with a Killing vector field, 
	{\it Trans. AMS.,} {\bf 349}, 1063-1080 (1997). 

%\vspace{-1.5mm}

\bibitem  {Sa-gre} S\'ANCHEZ M., 
Lorentzian manifolds admitting a Killing vector field, 
	{\it Nonlinear Anal.,} {\bf 30}, 643-654 (1997).

%\vspace{-1.5mm}

\bibitem  {Sa-ba} S\'ANCHEZ M., 
	Some remarks on Causality and Variational Methods in Lorentzian 	manifolds, 	
    {\it  Conf. Sem. Mat. Univ. Bari} {\bf 265} 1-12 (1997).

%\vspace{-1.5mm}



\bibitem{Se} SEIFERT H.J.,
Global connectivity by timelike geodesics, 
{\it Z. Naturefor.,} {\bf 22a}, 1356 (1970).

%%\vspace{-1.5mm}

\bibitem  {Sm} SMITH J.W., 
Lorentz structures on the plane, 
{\it Trans. AMS.} {\bf 95,} 226-237 (1960).

%\vspace{-1.5mm}

\bibitem{We} WEINSTEIN, T., {\it An introduction to Lorentz surfaces}, De Gruyter Exp. Math. {\bf 22} W. de G. N.Y. (1996).  

%\vspace{-1.5mm}

\bibitem{Wo}  WOLF J.A. {\it Spaces of constant curvature}, MacGraw-Hill, N.Y. (1967).






}

\end{thebibliography}
\end{document}